\documentclass[12pt]{article}
\usepackage{amssymb}
\usepackage{amsmath}
\DeclareSymbolFont{extraup}{U}{zavm}{m}{n}
\DeclareMathSymbol{\clubsuit}{\mathalpha}{extraup}{88}

\newcommand\F{{\cal F}}

\newtheorem{thm}{Theorem} \newtheorem{lem}{Lemma}

\title{Nearly Erd\H os-Ko-Rado theorems}
\author{ {\bf Gyula O.H. Katona}\\ R\'enyi Institute, HUN-REN\\
	Budapest Pf 127, 1364 Hungary \\ ohkatona@renyi.hu  \and {\bf Jian Wang}\\
Department of Mathematics, Taiyuan University of Technology\\
 Taiyuan, 030024, China\\ wangjian01@tyut.edu.cn}

\begin{document}
\date{}

\maketitle

\begin{abstract}
If a family $\mathcal{F}$ of $k$-element subsets of an $n$-element set is pairwise intersecting, $2k\leq n$ then $|\mathcal{F}|\leq {n-1\choose k-1}$ holds by the celebrated Erd\H{o}s-Ko-Rado theorem. But an intersecting family obviously satisfies the condition
$${\ell \choose 2}\leq \sum_{1\leq i<j\leq \ell}|F_i\cap F_j| $$
for any $\ell$ distinct members of the family. It has been proved in \cite{FKN} that even if ${\ell \choose 2}$ is replaced by
${\ell -1 \choose 2}+1$  the conclusion $|\mathcal{F}|\leq {n-1\choose k-1}$
remains valid for large $n$. However the 1 cannot be omitted, because there is a larger family satisfying that weaker condition. In the present paper we determine the largest size of the family under this weaker condition when $n$ is sufficiently large. All of these are treated in the more general setting of $t$-intersecting families.

{\it Key Words:} Extremal set theory, intersecting families,
Erd\H os-Ko-Rado theorem

\end{abstract}

\section{Introduction}
Let $[n]=\{ 1,2,\ldots , n \}  $ be an $n$-element set. ${[n] \choose k}$ denotes the family of
all $k$-element subsets of $[n]$. A family $\F \subset {[n] \choose k}$
 is called {\it intersecting} if any two members of $\F$ have a non-empty intersection.
Let us start with the
classic theorem of Erd\H{o}s, Ko and Rado.

\begin{thm} {\rm (\cite{EKR})} If $2k\leq n$ and $\F \subset {[n] \choose k}$ is intersecting then
$$|\F|\leq {n-1\choose k-1}.\eqno(1)$$
\end{thm}

Choose an integer $\ell \geq 2$ and take the following sum for $\ell$ distinct members of the family $\F$:
$$\sum_{1\leq i<j\leq \ell}|F_i\cap F_j|. $$
If ${\cal F}$ is intersecting then every term  is at least 1, therefore
$${\ell \choose 2}\leq \sum_{1\leq i<j\leq \ell}|F_i\cap F_j|$$
holds. Does this weaker condition imply (1)? Much more is true for large $n$.

\begin{thm} {\rm (\cite{FKN})} Let $2\leq k, \ell$ be integers and suppose that ${\cal F} \subset {[n]\choose k}$. If
$${\ell -1 \choose 2}+1\leq \sum_{1\leq i<j\leq \ell}|F_i\cap F_j|\eqno(2)$$
holds for any subfamily $\{ F_1, F_2, \ldots, F_{\ell}\}$ of distinct members of ${\cal F}$ and $n$ is sufficiently large then
$$|{\cal F}|\leq {n-1\choose k-1}.$$
\end{thm}
 Let us remark that the special case $\ell=3$ was proved by Kartal Nagy \cite{N} earlier.

This theorem is an obvious sharpening of the Erd\H os-Ko-Rado theorem,
However, the famous paper of Erd\H os, Ko and Rado \cite{EKR} contains a more general form of their theorem, but it holds only for large $n$.

A family  ${\cal F} \subset {[n]\choose k}$ is called $t$-{\it intersecting} ($1\leq t<k $)
if $|F_1\cap F_2|\geq t$ holds for every pair of members of  ${\cal F}$.

\begin{thm} {\rm (\cite{EKR})}
If  $ {\cal F} \subset {[n]\choose k}$ is $t$-intersecting and $n$ is sufficiently large  then
$$|{\cal F}|\leq {n-t\choose k-t}.$$
\end{thm}

Let us mention that Erd\H os, Ko and Rado established this statement for a rather large $n\geq n_0(k,t)$ whereas the exact bound $n_0(k,t)=(k-t+1)(t+1)$ is due to Frankl \cite{F2} for $t\geq 15$ and Wilson \cite{W} for smaller values of $t$.

An improvement of Theorem 3 was also found in \cite{FKN}.

If the family is $t$-intersecting then
$${\ell \choose 2}t\leq \sum_{1\leq i<j\leq \ell}|F_i\cap F_j|.$$

The left hand side here also can be decreased.

\begin{thm} {\rm (\cite{FKN})} Let $1 \leq t < k, 2 \leq \ell$  be integers and suppose that ${\cal F} \subset {[n]\choose k}$. If
	$${\ell \choose 2}(t-1)+{\ell -1 \choose 2}+1\leq \sum_{1\leq i<j\leq \ell}|F_i\cap F_j|\eqno(3)$$
	holds for any subfamily $\{ F_1, F_2, \ldots F_{\ell}\}$ of distinct members of ${\cal F}$ and $n$ is sufficiently large then
	$$|{\cal F}|\leq {n-t\choose k-t}.$$
\end{thm}

Our goal in the present paper is to investigate what happens if the bounds in (2) and (3) are further decreased by 1. Let us  first show the change in the special case when $t=1$.

\begin{thm}

Let $ 3\leq \ell, \ell-2\leq k$ be integers and suppose that ${\cal F} \subset {[n]\choose k}$. If
$${\ell -1 \choose 2}\leq \sum_{1\leq i<j\leq \ell}|F_i\cap F_j|$$
holds for any subfamily $\{ F_1, F_2, \ldots F_{\ell}\}$ of distinct members of ${\cal F}$ and $n$ is sufficiently large then
$$|{\cal F}|\leq {n-1\choose k-1}+{n-\ell +1\choose k-\ell +2}.$$
\end{thm}

And now let us see the general form.

\begin{thm}
	Let $1\leq t, 3\leq \ell , t+\ell -2\leq k$ be integers and suppose that ${\cal F} \subset {[n]\choose k}$. If
	$${\ell \choose 2}(t-1)+{\ell -1\choose 2}\leq \sum_{1\leq i<j\leq \ell}|F_i\cap F_j|\eqno(4)$$
	holds for any subfamily $\{ F_1, F_2, \ldots , F_{\ell}\}$ of distinct members of ${\cal F}$ and $n$ is sufficiently large then
	$$|{\cal F}|\leq {n-t\choose k-t}+{n-t-\ell +2 \choose k-t-\ell +3 }.\eqno(5)$$
\end{thm}

Let us note that these theorems are not valid for $\ell=2$.

In the proofs we will use the following theorem of Paul Erd\H os.

\begin{thm} {\rm (\cite{Er})} Let $\F\subset {[n]\choose k}$  and suppose that the family contains no $\ell$ paiwise disjoint members.  If $n$ is sufficiently large then
	$$|\F|\leq {n\choose k}-{n-\ell +1\choose k}.$$
\end{thm}

Let us note that in his paper Erd\H os posed a conjecture for the maximum size of the family for every $n$. It became known as the Erd\H os Matching Conjecture. It is still not fully solved, but there is  considerable progress, see \cite{F}, \cite{FK1} and \cite{FK2}.

Since Theorem 5 is a special case of Theorem 6, it is sufficient to prove the latter one.

\section{Proof of Theorem 6}

{\it Construction}. Take all the $k$-element subsets of $[n]$ containing the set $[t]$ and all the $k$-element subsets containing the set $[2, t+\ell -2]$, but not containing the element 1.
It is easy to see that their number is equal to the right hand side of (5).

Let us show that they satisfy condition (4). Choose $0\leq x\leq \ell$ sets from the subfamily containing $[t]$ and $\ell -x$ from the rest. Then the sum of the sizes of the pairwise intersections will be at least
$$f(x)={x \choose 2}t+x(\ell -x)(t-1)+{\ell -x\choose 2}(t+\ell -3).$$
This is a quadratic function on the real line with one minimum. Since $f(\ell -2)=f(\ell -1)$ these must be two minimum places if $x$ is restricted to be integer. Hence we have
$${\ell \choose 2}(t-1)+{\ell -1\choose 2}\leq f(x)$$
for integers.

{\it Upper bound.}
A {\it sunflower} is a family of subsets of the same size whose pairwise intersections give the same set. The common intersection is called the {\it kernel}. The disjoint rests are the {\it petals}. Let ${\cal S}(k,t,u)\subset {[n]\choose k}$ be a sunflower with a kernel $T$  of size $t$ and with $u$ petals.

{\it Case 1:} $\F$ contains a sunflower ${\cal S}(k,t,2k+\ell -2)$ as a subfamily.

This case will be settled by  a series of lemmas.
If ${\cal A} , {\cal B}\subset {[n]\choose k}$ we say that ${\cal A}$ and $ {\cal B}$ are {\it cross $r$-intersecting } $(r<k)$ if $A\in {\cal A}, B\in {\cal B}$ implies $r\leq |A\cap B|$.

\begin{lem} Let ${\cal A}, {\cal B}\subset {[n]\choose k}$ be $r$-intersecting ($r<k$) and $r+1$-cross intersecting families and $n$ is sufficiently large then
 $$	|{\cal A|+|\cal B}|\leq {n-r\choose k-r}\eqno(6)$$
 with equality only if one of the families is empty, otherwise both terms are $=O(n^{k-r-1}).$
 \end{lem}

 {\it Proof}. If one of the families is empty then (6) is obvious. Fix now one member of each of the families: $A\in {\cal A}, B\in {\cal B}$. Then any other member $C\in {\cal A}$ intersects $B$ in at least $r+1$ elements therefore the number of possible choices is at most ${k\choose r+1}{n-r-1\choose k-r-1}=O(n^{n-r-1})$. That is, $|{\cal A}|=O(n^{n-r-1})$. The same is true for $|{\cal B}|$, proving (6) and the rest of the statements.
  $\Box$

 \begin{lem} All members $F\in \F$ intersect the kernel $T$ of ${\cal S}(k,t,2k+\ell -2)$ in at least $t-1$ elements.
 \end{lem}
 	
 {\it Proof}. Since each $F\in \F$ can meet at most $k$  petals,  at least $k+\ell -2\geq \ell -1$ petals are disjoint to $F$. Let $P_1, P_2, \ldots , P_{\ell -1}$ be $\ell -1$ of them. Apply the condition of the theorem for $F$ and $P_1\cup T, P_2\cup T, \ldots , P_{\ell -1}\cup T$.
 The sum of the sizes of the pairwise intersections is
 $${\ell -1\choose 2}t+|F\cap T|(\ell -1).$$
		 This reaches the left hand side of (3) if $|F\cap T|\geq t-1. \ \Box$
		
 Now we conclude that  $\F$ can be decomposed into $t+1$ subfamilies: $\F (T)$ (the subfamily of members of $\F$ containing $T$), and the subfamilies $\F(T-\{ a\}, \bar{a})$ ( the subfamily of the members of $\F$ containing  $T-\{ a\}$ and not containing $a$), for every $a\in T$. Then
 $$|\F|=|\F (T)|+\sum_{a\in T}|\F(T-\{ a\}, \bar{a})|.\eqno(7)$$

 \begin{lem} $\F(T-\{ a\}, \bar{a})$ is $t-1+\ell -2$-intersecting.
 \end{lem}

  {\it Proof}.  Let $F, G\in \F(T-\{ a\}, \bar{a})$. These two sets can meet at most $2k$ petals of the sunflower, therefore there are petals $P_1, P_2, \ldots , P_{\ell -2}$ disjoint to them.
  Apply the condition of the theorem for the members $F, G, P_1, P_2, \ldots , P_{\ell -2}$.
  The sum of the sizes of their pairwise intersections is
  $${\ell -2\choose 2}t+2(t-1)(\ell -2)+|F\cap G|.\eqno(8)$$
 This reaches the left hand side of (3) if $|F\cap G|\geq t-1+\ell -2$.
 However, here $|F\cap G\cap T|=t-1$, this is why $|(F-T)\cap (G-T)|\geq \ell -2$ holds here. $ \Box$

 In other words, $\F(T-\{ a\}, \bar{a})-T$ is $\ell -2$-intersecting.

 \begin{lem} If $a, b\in T, a\not= b$ then  $\F(T-\{ a\}, \bar{a})-T$ and $\F(T-\{ b\}, \bar{b})-T$ are $\ell -1$-cross-intersecting.
 \end{lem}

 {\it Proof}.  Let $F\in \F(T-\{ a\}, \bar{a}), G\in \F(T-\{ b\}, \bar{b})-T$. These two sets can meet at most $2k$ petals of the sunflower, therefore there are petals $P_1, P_2, \ldots , P_{\ell -2}$ disjoint to them.
 Apply the condition of the theorem for the members $F, G, P_1, P_2, \ldots , P_{\ell -2}$.
 The sum of the sizes of their pairwise intersections is equal to (7), again. Hence we have
 $|F\cap G|\geq t-1+\ell -2$ as before. However, here  $|F\cap G \cap T|=t-2$, this is why
 $|(F-T)\cap (G-T)|\geq \ell -1$ holds here. $\Box$
 	
  \begin{lem} Either all subfamilies
  	 $\F(T-\{ a\}, \bar{a})\ (a\in T)$ are empty with at most one exception or
  	 $$\sum_{a\in T}|\F(T-\{ a\}, \bar{a})|\leq O(n^{k-t-\ell +2}).$$
 \end{lem}

 {\it Proof}. Suppose that $\F(T-\{ a\}, \bar{a})$ and $ \F(T-\{ b\}, \bar{b})$ are non-empty.
The members of $\F(T-\{ a\}, \bar{a})-T$ and $ \F(T-\{ b\}, \bar{b})-T$ are $k-t+1$-element subsets of an $n-t$-element set $[n]-T$. By Lemmas 1, 3, 4 and 5 both families  are $=O(n^{k-t-\ell +2})$, proving the lemma. $\Box$

$|\F (T)|\leq {n-t\choose k-t}$ is obvious. If $\F(T-\{ a\}, \bar{a})$ is non-empty only for one $a$ then we have
$$|\F(T-\{ a\}, \bar{a})|\leq {n-t-\ell +2 \choose k-t-\ell +3}$$
by Lemma 3.  Substituting these value into (7) we obtain the statement of the theorem in this case. On the other hand if there are more $a$'s  with non-empty
$\F(T-\{ a\}, \bar{a})$ then we obtain the upper bound
$$|\F|\leq  {n-t\choose k-t}+O(n^{k-t-\ell +2})$$
which is smaller for large $n$. Case 1 is settled.

{\it Case 2:} $\F$ contains no sunflower ${\cal S}(k,t,2k+\ell -2)$ as a subfamily.

\begin{lem} Let $T\subset [n]$ be a fixed $t$-element subset. If $\F$ contains no sunflower ${\cal S}(k,t,2k+\ell -2)$ with kernel $T$ then
	$$|\F (T)|=O(n^{k-t-1}).$$
\end{lem}

{\it Proof}. Our conditions imply that the family $\F (T) -T \subset{[n]-T\choose k-t}$ does not contain $2k+\ell -2$ pairwise disjoint members. By Theorem 7 we obtain
$$|\F (T)|=|\F (T)-T|\leq {n-t\choose k-t}-{n-t -(2k+\ell -2)+1\choose k-t}=O(n^{k-t-1}).$$
 $\Box$

Now we are ready to finish the proof also in Case 2. Choose the largest subfamily $G_1, G_2, \ldots , G_m\in \F$ satisfying $|G_i\cap G_j| \leq t-1\ (1\leq i<j\leq m).$
$m<\ell$ must hold, otherwise
$$\sum_{1\leq i<j\leq \ell}|G_i\cap G_j|\leq {\ell \choose 2}(t-1),$$
contradicting condition (4). Let $W=\cup_{i=1}^mG_i$. All members of $\F$ must intersect one $G_i$ in at least $t$ elements, therefore $|W\cap F|\geq t$ holds for every $F\in \F$. Consequently
$$\F = \bigcup_{T\subset W,\ |T|=t}\F(T)$$
holds.
Since $|W|\leq km\leq k(\ell -1)$, we have
$$|\F|\leq \sum_{T\subset W,\ |T|=t}|\F(T)|\leq {k(\ell -1) \choose t}O(n^{k-t-1})=O(n^{k-t-1})$$
by Lemma 6. But this  is smaller than the right hand side of (5) if $n$ is sufficiently large, finishing the proof also in this case.
 $\Box$

\section{Open problems}

{\it Problem 1.} Try to further decrease the lower bound in (3). For sake of simplicity let us consider here only the case $t=1$ and decrease the lower bound by $s$. Slightly modify the construction:
replace $[2, \ell -1]$ by $[2, \ell -s-1]$.
If $x$ sets containing the element 1 are chosen and $\ell -x$ from the other kind of sets then the sum of the sizes of the pairwise intersections is at least
$${x\choose 2}+{\ell -x\choose 2}(\ell -s-2).\eqno(9)$$
This is obviously larger than or equal to ${\ell -1 \choose 2}-s$ for $x=\ell, \ell -1 $ and $\ell -2$. However, for $x=\ell -3$ the following inequality is needed:
$${\ell -1 \choose 2}-s\leq {\ell -3 \choose 2}+3(\ell -s-2).$$
This is equivalent to $2s\leq \ell -1$. It means that the modified construction is good only when $2s\leq \ell -1$. But in this case the function (9) is the smallest for integer values at $\ell -2$.

On the other hand the proof of the upper bound in Theorem 6 also works for this case. This is why we can claim that the following theorem is true.

\begin{thm}
	Let $2s+1\leq \ell , \ell -s-1\leq k$ be integers and suppose that ${\cal F} \subset {[n]\choose k}$. If
	$${\ell -1\choose 2}-s\leq \sum_{1\leq i<j\leq \ell}|F_i\cap F_j|\eqno(10)$$
	holds for any subfamily $\{ F_1, F_2, \ldots , F_{\ell}\}$ of distinct members of ${\cal F}$ and $n$ is large enough then
	$$|{\cal F}|\leq {n-1\choose k-1}+{n-\ell +s+1 \choose k-\ell +s+2 }.$$
\end{thm}

Determine the largest family satisfying (10) with $\ell \leq 2s$ for any
subfamily $\{ F_1, F_2, \ldots , F_{\ell}\}$ of distinct members of ${\cal F}$ where $n$ is sufficiently large.

{\it Problem 2.} In our theorems it was supposed that $n$ is large. We did not even try to find a good bound for $n$ under which the statements are true. Find a good bound and determine the maximum of $|\F|$ for smaller $n$ satisfying the other conditions of our theorems.

{\it Problem 3.} In our theorems uniform families were considered that is $\F \subset {[n]\choose k}$ was supposed. Find analogous results for the non-uniform case when $\F \subset 2^{[n]}$. The case $\ell =3$ has been solved by Kartal Nagy in \cite{N}.

\vspace{3pt}
{\noindent \bf Acknowledgement.} The second author was supported by
National Natural Science Foundation of China Grant no. 12471316.

\end{document}